\documentclass{amsart}
\usepackage{setspace}

\usepackage{amsmath}
\usepackage{amsthm}
\usepackage{amssymb}

\usepackage[latin1]{inputenc} 

\newtheorem{theorem}{Theorem}[section]

\newtheorem{proposition}[theorem]{Proposition}

\numberwithin{equation}{section}

\begin{document}
\title[Some generalizations of the functions $\tau$ and $\tau^{\left(e\right)}$ in algebraic number fields]
{Some generalizations of the functions $\tau$ and $\tau^{\left(e\right)}$ in algebraic number fields}

\author{Nicu\c{s}or Minculete}
\address{Faculty of Mathematics and Computer Science, Transilvania University\\
 Iuliu Maniu street 50, Bra\c{s}ov 500091, Romania}
\email{minculete.nicusor@unitbv.ro}

\author{Diana Savin}
\address{Faculty of Mathematics and Computer Science, Ovidius University\\
Bd. Mamaia 124, 900527, Constan\c{t}a, Romania}
\email{savin.diana@univ-ovidius.ro; dianet72@yahoo.com}

\date{}
\begin{abstract}
In this paper, we generalize the arithmetic functions $\tau$ and $\tau^{\left(e\right)}$ in algebraic number fields and we find some properties of these functions.
\end{abstract}
\maketitle
 
\section{ Preliminaries}
\noindent Let $n\geq 1$ be a natural number and we define by $\tau\left(n\right)$ the number of divisors of $n$. We note that $\tau\left(1\right)=1$; if $p$ is a prime number, then $\tau\left(p\right)=2$ and $\tau\left(p^a\right)=a+1$. It is known that the function $\tau$ is multiplicative, but not completely multiplicative and for $n=p_1^{a_1}p_2^{a_2}\cdot\cdot\cdot p_r^{a_r}$, $n>1$, we have the relation $\tau\left(n\right)=\left(a_1+1\right)\left(a_2+1\right)\cdot\cdot\cdot\left(a_r+1\right).$\\
\indent A natural number $d$ is a unitary divisor of a number $n$ if $d$ is a divisor of $n$ and if $d$ and $\frac{n}{d}$ are coprime, so they have no common factor other than $1.$ If $\tau^*\left(n\right)$ is the number of the unitary divisors of $n$, then for $n=p_1^{a_1}p_2^{a_2}\cdot\cdot\cdot p_r^{a_r}$, $n>1$, we obtain the relation $\tau^*\left(n\right)=2^r=2^{\omega\left(n\right)},$ where $\omega\left(n\right)$ denote the number of distinct prime factors of $n$.\\
\indent The notion of \textit{exponential divisor} was introduced by Subbarao \cite{Subbarao} in the following way: if we consider a natural number $n>1$ which can be written in canonical form as $n=p_1^{a_1}p_2^{a_2}\cdot\cdot\cdot p_r^{a_r}$, the natural number $d=\prod_{i=1}^{r}p_i^{b_i}$ is called \textit{exponential divisor} or \textit{e-divisor} of $n=\prod_{i=1}^{r}p_i^{a_i}>1$ if $b_i|a_i$, for every $i\in \{1,...,r\}.$ We write $d|_{\left( e\right)}n$.\\
\indent Let $\tau^{\left( e\right)}\left(n\right)$ be the number of exponential divisors of $n$. By convention $\tau^{\left( e\right)}\left(1\right)=1.$ It is easy to see that the function $\tau^{\left(e\right)}$ is multiplicative and for $n=p_1^{a_1}p_2^{a_2}\cdot\cdot\cdot p_r^{a_r}>1$ we have $\tau^{\left( e\right)}\left(n\right)=\tau\left( a_1\right)\tau\left( a_2\right)\cdot\cdot\cdot\tau\left( a_r\right).$ 
It is obvious that $\tau^{\left( e\right)}\left(n\right)\leq\tau\left(n\right)$, for all $n\in\mathbb{N}^{*}.$\\
\indent Other properties of the number of the exponential divisors of $n$ can be found in the papers \cite{Min_2}, \cite{Sandor_2} and \cite{Toth}.\\
\indent Let $n\in$$\mathbb{N}$, $n\geq2$ and let $K$ be an algebraic number field of degree $[K :
\mathbb{Q}] = n.$ Let $p$ be a prime positive integer. Let $\mathcal{O}_{K}$ be the ring of integers of the field $K$. It is known that $\mathcal{O}_{K}$ is a Dedeking ring. According to the fundamental theorem about Dedekind rings, there exist positive integers $g$ and $e_{i},$ $i = \overline{1,g}$ and the different ideals $P_{1},$ $P_{2},$..., $P_{g}$$\in$Spec$\left(\mathcal{O}_{K}\right)$  such that
$$p\mathcal{O}_{K}=P^{e_{1}}_{1}\cdot P^{e_{2}}_{2}\cdot...\cdot P^{e_{g}}_{g}.$$
This decomposition is unique, except the order of factors.
The integer $e_{i}$ ($i = \overline{1,g}$) is called \textit{the ramification index} of $p$ at the ideal $P_{i}$ and the
degree $f_{i}$  of the following extension of fields $\left(\mathcal{O}_{K}/P_{i}\right)/\left( \mathbb{Z}/p\mathbb{Z}\right)$ i. e. $f_{i}=\left[\mathcal{O}_{K}/P_{i} : \mathbb{Z}/p\mathbb{Z}\right]$ is called \textit{the residual degree} of $p$ (see \cite{ireland}). There are the following formulas (see \cite{ireland}):

\begin{proposition}
a) With the above notations, we have
$$\sum_{i=1}^{g}e_{i}f_{i}=[K :\mathbb{Q}] = n$$
b) If $\mathbb{Q} \subset K$ is a Galois extension, then $e_{1}=e_{2}=...=e_{g}$ (say to $e$), $f_{1}=f_{2}=...=f_{g}$ (say to $f$).
So, $efg=n.$
\end{proposition}

We recall now how a prime ideal of $\mathbb{Z}$ is decomposing into the ring of integers of a quadratic field, into the ring of integers of a cyclotomic field, respectively, into the ring of integers of a Kummer field, the results we will use in the next section.

\begin{theorem}[\cite{ireland}]
Let $d\neq 0,1$   {be a square free   integer.}   {Let} $\mathcal{O}_{K}$   {be the ring of integers of the quadratic field} $K=\mathbb{Q}(\sqrt{d})$   {and} $\Delta_{K}$   {be the discriminant of} $K.$   {Let} $p$   {be an odd prime integer. Then, we have}:
\begin{enumerate}
\item $p$   {is ramified in} $\mathcal{O}_{K}$ if and only if
$p  | \Delta_{K}$. In this case 
$p\mathcal{O}_{K}=(p   , \sqrt{d}   )^2$;
\item $p$   splits totally   in $\mathcal{O}_{K}$  if and only if the Legendre symbol
  $(\frac{\Delta_{K}}{p})=1$. In this case $p\mathcal{O}_{K}=P_1 \cdot  {P_2},$   {where} $P_1$ and $ P_2$   {are distinct prime ideals in} $\mathcal{O}_{K};$
\item $p$   {is inert in} $\mathcal{O}_{K}$ {if and only if} the Legendre symbol
 $(\frac{\Delta_{K}}{p})=-1$; 
\item 
the prime 2 is ramified in $\mathcal{O}_{K}$  {if and only if}
 $d\equiv 2$ (mod $4$) or $d\equiv 3$ (mod $4$) 
  In the first case    $2\mathcal{O}_{K}=(2,\sqrt{d} )^{2}$, while in the second case
 $2\mathcal{O}_{K}=(2,1+\sqrt{d} )^{2};$
\item   
the prime $2$  splits totally   in $\mathcal{O}_{K}$  {if and only if} 
$d\equiv 1$ (mod 8). In this case     $2\mathcal{O}_{K}=P_1 \cdot P_2$,   {where} $P_1, P_2$
  {are distinct prime ideals in} $\mathcal{O}_{K}$, with $P_1=(2, \frac{1+\sqrt{d}}{2})$;
\item 
the prime $2 $   {is inert  in} $\mathcal{O}_{K}$  {if and only if} 
$d\equiv 5$ (mod $8$). 
\end{enumerate}
\end{theorem}

Let $n$ be a positive integer, $n\geq2$ and let $\xi$  be a primitive root of unity of degree $n.$ It is known that the ring of integers of the cyclotomic field $\mathbb{Q}\left(\xi\right)$ is $\mathbb{Z}\left[\xi\right].$\\
We recall now a result of Dedekind.
\begin{theorem} [\cite{ireland}, \cite{Savin}]
Let K be an algebraic number field and p a natural prime number. Then p is ramified in $\mathcal{O}_{K}$ if and only if $p|\Delta_K.$
\end{theorem}
\begin{theorem}[\cite{Hilbert}, \cite{ireland}]
Let $n$ be an integer, $n \geq 3,$ $\xi$  be a primitive root of unity of order
$n$ and let the cyclotomic field $\mathbb{Q}\left(\xi\right)$. If $p$ is a prime integer
which does not divide $n$ and $f$ is the smallest positive integer such that $p^{f}$$\equiv$$1$ (mod $n$),
then
$$p\mathbb{Z}\left[\xi\right]=P_{1}\cdot P_{2}\cdot...\cdot P_{g},$$
where $g=\frac{\varphi\left(n\right)}{f}$ and $P_{j}$, $j =\overline{1,g}$ are different prime ideals in the ring $\mathbb{Z}\left[\xi\right]$ and $\varphi$ is the Euler's function.

\end{theorem}
\begin{theorem}[\cite{Hilbert}, \cite{ireland}]
Let $p$ be an odd prime positive integer and let $\xi$ be a primitive root of unity of order $p$. Then the ideal $p\mathbb{Z}\left[\xi\right]$ is ramified totally in $\mathbb{Z}\left[\xi\right]$.
\end{theorem}

\begin{theorem}[\cite{Hilbert}] \textit{Let }$\xi $\textit{\ be a primitive
root of the unity of }$l-$\textit{order, where }$l$\textit{\ is a prime
natural number and let} $\mathcal{O}_{L}.$ \textit{be the ring of integers of the Kummer
field} $\mathbb{Q}(\xi ,\sqrt[l]{\mu })$ . \textit{A prime ideal} $P$\textit{\ in the
ring} $\mathbb{Z}[\xi ]$\textit{\ is in }$\mathcal{O}_{L}$\textit{\ in one of the
situations:}\\
\textit{i) It is equal with the }$l-$\textit{power of a prime ideal from} $%
\mathcal{O}_{L}, $ \textit{if the} $l-$\textit{power character }$\left( \frac{\mu }{P}%
\right) _{l}=0;$\\
\textit{ii) It is a prime ideal in} $\mathcal{O}_{L}$, \textit{if} $\left( \frac{\mu }{P}%
\right) _{l}=$\textit{\ a root of order }$l$\textit{\ of unity, different
from }$1$.\\
\textit{iii) It decomposes in }$l$\textit{\ different prime ideals from} $\mathcal{O}_{L}$%
\textit{if} $\left( \frac{\mu }{P}\right) _{l}=1.$ 
\end{theorem}

\section{The generalization of the arithmetic functions $\tau$ and $\tau^{\left(e\right)}$ in algebraic number fields}
\label{sectionarithmetic}

Let $K$ be an algebraic number field and let $\mathcal{O}_{K}$ be the ring of integers of $K$. We denote with $\mathbb{J}$ the set of ideals of the ring $\mathcal{O}_{K}.$ Let $I$ be an ideal of the ring $\mathcal{O}_{K}.$ Since $\mathcal{O}_{K}$ is a Dedekind ring, $\left(\exists!\right)$ $g$$\in$$\mathbb{N}^{*},$ the different ideals $P_{1},$ $P_{2},$..., $P_{g}$$\in$Spec$\left(\mathcal{O}_{K}\right)$ and the positive integer numbers $e_{1},$ $e_{2},$...,$e_{g}$$\in$$\mathbb{N}^{*}$ such that $I=P^{e_{1}}_{1}\cdot P^{e_{2}}_{2}\cdot...\cdot P^{e_{g}}_{g}.$ A \textit{divisor} of the ideal $I$ is of the following form $d_{I}=P^{a_{1}}_{1}\cdot P^{a_{2}}_{2}\cdot...\cdot P^{a_{g}}_{g},$
where $a_{1},$ $a_{2},$...,$a_{g}$$\in$$\mathbb{N},$ with $a_{i}$$\leq$$e_{i},$ for $\left(\forall\right)$ $i=\overline{1,g}.$ An \textit{exponential divisor} of the ideal $I$ is of the following form
$d^{\left(e\right)}_{I}=P^{b_{1}}_{1}\cdot P^{b_{2}}_{2}\cdot...\cdot P^{b_{g}}_{g},$
where $b_{1},$ $b_{2},$...,$b_{g}$$\in$$\mathbb{N}^{*},$ with $b_{i}$ $|$ $e_{i},$ for $\left(\forall\right)$ $i=\overline{1,g}.$

We extend the functions $\tau$ and $\tau^{\left(e\right)}$ to ideals in the following way:
$\tau : \mathbb{J}\rightarrow \mathbb{C},\tau\left(I\right)$ is the number of divisors of the ideal $I$ and $\tau^{\left(e\right)} : \mathbb{J}\rightarrow\mathbb{C}$, $\tau^{\left(e\right)}\left(I\right)$ is the number of exponential divisors of the ideal $I$. So, we have $$\tau^{\left(e\right)}\left(I\right)=\sum_{d^{\left(e\right)}_{I}|_{\left(e\right)}I}1.$$

We quickly obtain the following formulas:
\begin{proposition}
With the above notations, we have
$$\tau\left(I\right)=\left(e_{1}+1\right)\cdot\left(e_{2}+1\right)\cdot...\cdot\left(e_{g}+1\right) $$ and $$\tau^{\left(e\right)}\left(I\right)=\tau\left(e_{1}\right)\cdot \tau\left(e_{2}\right)\cdot...\cdot \tau\left(e_{g}\right).$$
\end{proposition}

If $p$ is a prime positive integer, then we study the number of the exponential divisors of the ideal $p\mathcal{O}_{K}$ and and their shape.\\
Let $\mathbb{Q} \subset K$ be a Galois extension of degree $\left[K:\mathbb{Q}\right]=n<\infty$, let $p$ be a prime positive integer and let $\mathcal{O}_{K}$ be the ring of integers of the field $K$. If $p\mathcal{O}_{K}$ has decomposition into product of prime ideals, given by $p\mathcal{O}_{K}=P^{e_{1}}_{1}\cdot P^{e_{2}}_{2}\cdot...\cdot P^{e_{g}}_{g}$, where $e_{i}\in \mathbb{N}^*,$ $i = \overline{1,g}$ and the different ideals $P_{1},$ $P_{2},$..., $P_{g}$$\in$Spec$\left(\mathcal{O}_{K}\right),$  then from the fact that $\mathbb{Q} \subset K$ is a Galois extension we obtain $e_1=e_2=...=e_g$. Therefore, we deduce $p\mathcal{O}_{K}=P^{e_{1}}_{1}\cdot P^{e_{1}}_{2}\cdot...\cdot P^{e_{1}}_{g}$, so we find the relation
$$\tau^{\left( e\right)}\left(p\mathcal{O}_{K}\right)=\left(\tau\left(e_{1}\right)\right)^g.$$
\indent We consider $K$ an algebraic number field of degree $[K :
\mathbb{Q}] = q$, where $q$ is a prime positive integer. If $\mathbb{Q} \subset K$ is a Galois extension and $\mathcal{O}_{K}$ is the ring of integers of the field $K$, then $$\tau^{\left( e\right)}\left(p\mathcal{O}_{K}\right)\in\{1, 2\},$$ because applying Propositon 1.1 (b) we distinguish the following cases:
i)  $p\mathcal{O}_{K}\in\text{Spec}\left(\mathcal{O}_{K}\right)$, implies that $\tau^{\left( e\right)}\left(p\mathcal{O}_{K}\right)=1;$
ii) $p\mathcal{O}_{K}=P_1\cdot P_2\cdot \cdot\cdot P_q$, where $P_i\in\text{Spec}\left(\mathcal{O}_{K}\right)$, for all $i = \overline{1,q}$, prove that $\tau^{\left( e\right)}\left(p\mathcal{O}_{K}\right)=1$ and iii) $p\mathcal{O}_{K}=P^q$, with $P\in\text{Spec}\left(\mathcal{O}_{K}\right)$, implies that $\tau^{\left( e\right)}\left(p\mathcal{O}_{K}\right)=2,$ because the ideal $p\mathcal{O}_{K}$ has only exponential divisors $P$ and $P^q$.\\
\indent S\'andor (\cite{Sandor}) established the following inequality:
$$2^{\omega\left(n\right)}\leq \tau^{\left( e\right)}\left(n\right)\leq2^{\Omega\left(n\right)},$$ for all natural numbers $n\geq 2$, such that $n$ is a perfect square, where $\omega\left(n\right)$ represents the number of distinct prime factors of $n$ and $\Omega\left(n\right)$ represents the number of distinct prime factors of $n$ together with their multiplicities. We study this inequality for ideals in the rings of integers of an algebraic number field. 
\begin{proposition}
Let $K$ be an algebraic number field of degree $[K :
\mathbb{Q}] = n<\infty$ and its ring of integers $\mathcal{O}_{K}$. Then, for every ideal I of the ring $\mathcal{O}_{K}$, we have
$$2^{\omega\left(I^2\right)}\leq \tau^{\left( e\right)}\left(I^2\right)<2^{\Omega\left(I^2\right)},$$
where $\omega\left(I^2\right)$ represents the number of distinct prime factors of $I^2$ ($\omega\left(I^2\right)=\omega\left(I\right)$) and $\Omega\left(I^2\right)$ represents the number of distinct prime factors of $I^2$ together with their multiplicities. 
\end{proposition}
\begin{proof}
Since the ring of integers $\mathcal{O}_{K}$ is a Dedekind ring and $I$ is a ideal of $\mathcal{O}_{K}$ there are unique $P_{1}, P_{2},..., P_{g}\in\text{Spec}\left(\mathcal{O}_{K}\right)$ and $a_1, a_2,...,a_g\in\mathbb{N}^*$ such that $I=P^{a_{1}}_{1}\cdot P^{a_{1}}_{2}\cdot...\cdot P^{a_{g}}_{g}$. It follows that $I^2=P^{2a_{1}}_{1}\cdot P^{2a_{1}}_{2}\cdot...\cdot P^{2a_{g}}_{g}$, so $\tau^{\left(e\right)}\left(I^2\right)=\tau\left(2a_{1}\right)\cdot \tau\left(2a_{2}\right)\cdot...\cdot \tau\left(2a_{g}\right).$ But, we have $\tau\left(2a_{i}\right)\geq 2,$ for all $i = \overline{1,g}$. Therefore we find $\tau^{\left( e\right)}\left(I^2\right)\geq 2^g$ and because $\omega\left(I^2\right)=g$, we deduce the inequality $\tau^{\left( e\right)}\left(I^2\right)\geq 2^{\omega\left(I^2\right)}.$
By mathematical induction, it easy to see that $\tau\left(n\right)<2^{2n-1}$, for every $n\geq 1$. Since $\tau\left(mn\right)\leq\tau\left(m\right)\tau\left(n\right)$, for every $m,n\in\mathbb{N}^*$, we deduce that $\tau\left(2a_i\right)\leq\tau\left(2\right)\tau\left(a_i\right)=2\tau\left(a_i\right)<2^{2a_i}$, for $a_i\geq 1$ and $i = \overline{1,g}$. This means that $\tau^{\left( e\right)}\left(I^2\right)=\tau\left(2a_{1}\right)\cdot \tau\left(2a_{2}\right)\cdot...\cdot \tau\left(2a_{g}\right)< 2^{2a_1+2a_2+...+2a_g}=2^{\Omega\left(I^2\right)}.$ Therefore, the statement is true.
\end{proof}
\begin{proposition}
For every number $n\in\mathbb{N}^*$, we have $$\tau\left(n\right)\leq 2^{n-1}.$$ 
\end{proposition}
\begin{proof} We consider the proposition $P\left(n\right):\tau\left(n\right)\leq 2^{n-1},$ for all $n\geq 1$.
Verify that $P\left(1\right)$ is true: $\tau\left(1\right)=1\leq 1=2^{1-1}.$ We assume that $P\left(k\right)$ is true, for all $k\leq n-1$. We prove that the proposition $P\left(n\right)$ is true. We identify two cases: i) If $n=q$, where $q$ is a prime number, then we have $\tau\left(n\right)=2\leq 2^{n-1}$; ii) If $n>2$ is a composite number, then we deduce that there are the numbers $d_1,d_2\in\mathbb{N}$, $2\leq d_1, d_2\leq n-1$, $\left(d_1, d_2\right)=1$, such that $n=d_1d_2$. Since $d_1, d_2\leq n-1$, we get $\tau\left(d_1\right)\leq 2^{d_1-1}$, $\tau\left(d_2\right)\leq 2^{d_2-1}$. Because the function $\tau$ is multiplicative, then we obtain $\tau\left(n\right)=\tau\left(d_1\right)\cdot\tau\left(d_2\right)\leq 2^{d_1-1}\cdot 2^{d_2-1}=2^{d_1+d_2-2}\leq 2^{d_1\cdot d_2-1}=2^{n-1}$, so $P\left(n\right)$ is true. According to the principle of mathematical induction, we have that $P\left(n\right)$ is true, for every $n\geq 1$.
\end{proof}
Next, we obtain a refinement of Proposition 2.2, given by the following:
\begin{proposition}
Let $K$ be an algebraic number field of degree $[K :
\mathbb{Q}] = n<\infty$ and its ring of integers $\mathcal{O}_{K}$. Then, for every ideal I of the ring $\mathcal{O}_{K}$, we have
$$2^{\omega\left(I^2\right)}\leq \tau^{\left( e\right)}\left(I^2\right)\leq2^{\Omega\left(I^2\right)-\omega\left(I^2\right)}.$$ 
\end{proposition} 
\begin{proof}Similarly, as in the proof from Proposition 2.2 and using the inequality from Proposition 2.3, we find the statement.
\end{proof}
We also obtain the following inequality involving the function $\tau$, the inequality that we will apply later.
\begin{proposition}
For every number $m,n\in\mathbb{N}, m\geq 2$, we have $$ \left(m+1\right)^n+1\geq\left[\tau\left(m\right)\right]^n+ 2^{n}.$$ 
\end{proposition}
\begin{proof} We consider the proposition $P\left(n\right):\left(m+1\right)^n+1\geq[\tau\left[\left(m\right)\right]^n+ 2^{n},$ for all $n\in\mathbb{N}$.
It is obviously that $P\left(0\right)$ is true. We check that $P\left(1\right)$ is true: the inequality of the statement becomes $m\geq\tau\left(m\right)$, which is true, for all $m\geq 2$. We assume that $P\left(n\right)$ is true, i. e. $\left(m+1\right)^n+1\geq[\tau\left[\left(m\right)\right]^n+ 2^{n}.$ Since, we have $\left(m+1\right)^{n+1}+1=\left(m+1\right)\left(m+1\right)^{n}+1\geq \left(m+1\right)\cdot\{\left[\tau\left(m\right)\right]^n+ 2^{n}-1\}+1=m\cdot\left[\tau\left(m\right)\right]^n+\left[\tau\left(m\right)\right]^n+2^{n}+m\left(2^{n}-1\right)\geq\tau\left(m\right)\cdot\left[\tau\left(m\right)\right]^n+2^{n}+2^{n}=\left[\tau\left(m\right)\right]^{n+1}+ 2^{n+1}$, because $m\left(2^{n}-1\right)\geq 0$, for every $n\geq 1$ and $m\geq\tau\left( m\right)\geq2$, for all $m\geq 2$. According to the principle of mathematical induction, we have that $P\left(n\right)$ is true, for every $n\in\mathbb{N}$.
\end{proof}
In \cite{Minculete}, the first author gave the following result:
\begin{proposition}
For any $n\geq 1$, the following inequality $$\tau\left(n\right)+1\geq\tau^{\left(e\right)}\left(n\right)+\tau^*\left(n\right)$$ 
occurs.
\end{proposition}
Next, we extend the above result to the followings:
\begin{proposition}
Let $K$ be an algebraic number field of degree $[K :
\mathbb{Q}] = n<\infty$ such that $\mathbb{Q} \subset K$ is a Galois extension. Let $\Delta_K$ be the discriminant of the field $K.$ If $p$ is a prime positive integer such that $p\not{|}\Delta_K$, then $$\tau\left(p\mathcal{O}_{K}\right)+1=\tau^{\left( e\right)}\left(p\mathcal{O}_{K}\right)+2^{\omega\left(p\mathcal{O}_{K}\right)}.$$ 
\end{proposition}
\begin{proof}
Since $p\not{|}\Delta_K$ it follows that $p$ is unramified in $\mathcal{O}_{K}$ (according to Theorem 1.3). Therefore there are unique $P_{1}, P_{2},..., P_{g}\in\text{Spec}\left(\mathcal{O}_{K}\right)$, $P_i\neq P_j$, for all $i,j = \overline{1,g}, i\neq j$, such that $p\mathcal{O}_{K}=P_{1}\cdot P_{2}\cdot...\cdot P_{g}$. This means that $\tau\left(p\mathcal{O}_{K}\right)=\left(1+1\right)\cdot\left(1+1\right)\cdot...\cdot\left(1+1\right)=2^g,  2^{\omega\left(p\mathcal{O}_{K}\right)}=2^g$ and $\tau^{\left( e\right)}\left(p\mathcal{O}_{K}\right)=1$. Therefore, we deduce the statement. 
\end{proof}
\begin{proposition}
Let $K$ be quadratic field. If $p$ is a prime positive integer, then $$\tau\left(p\mathcal{O}_{K}\right)+1=\tau^{\left( e\right)}\left(p\mathcal{O}_{K}\right)+2^{\omega\left(p\mathcal{O}_{K}\right)}.$$ 
\end{proposition}
\begin{proof}
According Theorem 1.2, we have three cases:
i) $p\mathcal{O}_{K}=P_{1}\cdot P_{2},$ where $P_{1}, P_{2}$$\in\text{Spec}\left(\mathcal{O}_{K}\right)$. This means that $\tau\left(p\mathcal{O}_{K}\right)=4,  2^{\omega\left(p\mathcal{O}_{K}\right)}=4$ and $\tau^{\left( e\right)}\left(p\mathcal{O}_{K}\right)=1$. Therefore, the inequality of the statement is true;\\
ii) $p\mathcal{O}_{K}=P_{1}^2$, where $P_{1}$$\in\text{Spec}\left(\mathcal{O}_{K}\right).$ It follows that $p$ is ramified totally in $\mathcal{O}_{K}.$ This implies $\tau\left(p\mathcal{O}_{K}\right)=3, 2^{\omega\left(p\mathcal{O}_{K}\right)}=2$ and $\tau^{\left( e\right)}\left(p\mathcal{O}_{K}\right)=2$. Therefore, the equality of the statement is true;\\
iii) $p\mathcal{O}_{K}=P\in\text{Spec}\left(\mathcal{O}_{K}\right)$, which proves that $p$ is inert in $\mathcal{O}_{K}$. Consequently, we obtain $\tau\left(p\mathcal{O}_{K}\right)=2, 2^{\omega\left(p\mathcal{O}_{K}\right)}=2$ and $\tau^{\left( e\right)}\left(p\mathcal{O}_{K}\right)=\tau\left(1\right)=1$. Therefore, the equality of the statement is true.
\end{proof}
Now we study what happens when $p|\Delta_K.$ We obtain the following result. 
\begin{proposition}
Let $K$ be an algebraic number field of degree $[K :
\mathbb{Q}] = n<\infty,$ $n\geq3$ such that $\mathbb{Q} \subset K$ is a Galois extension. Let $\Delta_K$ be the discriminant of the field $K.$ If $p$ is a prime positive integer such that $p|\Delta_K,$ then $$\tau\left(p\mathcal{O}_{K}\right)+1 \geq \tau^{\left( e\right)}\left(p\mathcal{O}_{K}\right)+2^{\omega\left(p\mathcal{O}_{K}\right)}.$$ 
\end{proposition}
\begin{proof}
Since $p|\Delta_K$ it follows that $p$ is ramified in $\mathcal{O}_{K}$ (according to Theorem 1.3). Therefore there are unique $P_{1}, P_{2},..., P_{g}\in\text{Spec}\left(\mathcal{O}_{K}\right)$, $P_i\neq P_j$, for all $i,j = \overline{1,g}, i\neq j$ and there is $e$$\in$$\mathbb{N},$ $e\geq2$ such that $p\mathcal{O}_{K}=P^{e}_{1}\cdot P^{e}_{2}\cdot...\cdot P^{e}_{q}$ and $efg=n$. This implies that $\tau\left(p\mathcal{O}_{K}\right)=\left(e+1\right)^g, 2^{\omega\left(p\mathcal{O}_{K}\right)}=2^g$ and $\tau^{\left( e\right)}\left(p\mathcal{O}_{K}\right)=\left(\tau\left(e\right)\right)^g$. Applying Proposition 2.5 we obtain:
$$\left(e+1\right)^g +1\geq \left(\tau\left(e\right)\right)^g + 2^g.$$
This inequality is equivalent with
$$\tau\left(p\mathcal{O}_{K}\right)+1 \geq \tau^{\left( e\right)}\left(p\mathcal{O}_{K}\right)+2^{\omega\left(p\mathcal{O}_{K}\right)}.$$ 
\end{proof}

If $p$ is totally ramified in $\mathcal{O}_{K}$, then we obtain the strict inequality, which is true in any extension of fields, not only in the Galois extension with the degree $\geq3$.
\begin{proposition}
Let $K$ be an algebraic number field of degree $[K :
\mathbb{Q}] = n\geq 3$. If p is a prime positive integer such that $p$ is totally ramified in $\mathcal{O}_{K}$, then we have the following:
$$\tau\left(p\mathcal{O}_{K}\right)+1>\tau^{\left( e\right)}\left(p\mathcal{O}_{K}\right)+2^{\omega\left(p\mathcal{O}_{K}\right)}.$$ 
\end{proposition}
\begin{proof}
Since $p$ is ramified totally in $\mathcal{O}_{K}$, we get $p\mathcal{O}_{K}=P^g$, where $P\in\text{Spec}\left(\mathcal{O}_{K}\right)$. If $f$ is the residual degree, then using Proposition 1.1, we find $e=n, g=1, f=1$, which implies  $p\mathcal{O}_{K}=P$, with $P\in\text{Spec}\left(\mathcal{O}_{K}\right)$. Consequently, we obtain $\tau\left(p\mathcal{O}_{K}\right)=n+1$, $2^{\omega\left(p\mathcal{O}_{K}\right)}=2$ and $\tau^{\left( e\right)}\left(p\mathcal{O}_{K}\right)=\tau\left(n\right)$. It is easy to see that $\tau\left(n\right)<n$, for $n\geq 3$, this shows that the inequality of the statement is true. 
\end{proof}

In \cite{Minculete}, N. Minculete (the first author) obtained the following result:

\begin{theorem}
For each integer number $n=p^{a_{1}}_{1}p^{a_{2}}_{2}...p^{a_{r}}_{r}$$>$$1$ (where $p_{1}, p_{2},...,p_{r}$ are different prime positive integers and $a_{1}, a_{2},...,a_{r}$$\in$$\mathbb{N}^{*}$), we have the following inequality:
$$\tau\left(n\right)\geq \tau^{\left(e\right)}\left(n\right) + \frac{\tau\left(n\right)}{\omega\left(n\right)}\cdot\left(\frac{1}{a_{1}+1}+\frac{1}{a_{2}+1}+...+\frac{1}{a_{r}+1}\right).$$
\end{theorem}

Now, we generalize this result for ideals in the ring of integers of an algebraic number field $K,$ when $\mathbb{Q} \subset K$ is a Galois extension. We obtain:

\begin{theorem}  
Let $n$$\in$$\mathbb{N}$, $n\geq 2$ and let $K$ be an algebraic number field of degree $[K :
\mathbb{Q}] = n$ such that $\mathbb{Q} \subset K$ is a Galois extension. Let $p$ be a prime positive integer and let the ideal
$p\mathcal{O}_{K}=P^{e}_{1}\cdot P^{e}_{2}\cdot...\cdot P^{e}_{g},$ where $g$$\in$$\mathbb{N}^{*},$ $e$$\in$$\mathbb{N}^{*},$ and $P_{1},$ $P_{2},$..., $P_{g}$ are prime different ideals in the ring of integers $\mathcal{O}_{K}.$
Then, we have the following inequality:
\begin{equation}{\label{3.1}}
\tau\left(p\mathcal{O}_{K}\right)\geq \tau^{\left(e\right)}\left(p\mathcal{O}_{K}\right) + \frac{\tau\left(p\mathcal{O}_{K}\right)}{\omega\left(p\mathcal{O}_{K}\right)}\cdot \frac{g}{e+1}.
\end{equation}
Moreover, the equality occurs if and only if $p\mathcal{O}_{K}$ is a prime ideal or $p\mathcal{O}_{K}=P^{2},$ with $P$$\in$Spec$\left(\mathcal{O}_{K}\right).$ 
\end{theorem}
\begin{proof}
We have: $\tau\left(p\mathcal{O}_{K}\right)=\left(e+1\right)^{g},$ $\tau^{\left(e\right)}\left(p\mathcal{O}_{K}\right)=\left(\tau\left(e\right)\right)^{g},$ $\omega\left(p\mathcal{O}_{K}\right)=g.$
The inequality (2.1) is equivalent with 
$$\left(e+1\right)^{g}\overset{?}{\geq} \left(\tau\left(e\right)\right)^{g} + \left(e+1\right)^{g-1} \Leftrightarrow$$
\begin{equation}{\label{3.2}}
\left(e+1\right)^{g-1}\cdot e\overset{?}{\geq} \left(\tau\left(e\right)\right)^{g}. 
\end{equation}
In order to prove the inequality (2.2), we consider three cases:\\
\textbf{Case 1.} If $p\mathcal{O}_{K}$$\in$Spec$\left(\mathcal{O}_{K}\right),$ then $e=g=1$ and we obtain the equality in (2.2) and also in (2.1).\\
\textbf{Case 2.} If $p\mathcal{O}_{K}=P^{2},$ with $P$$\in$Spec$\left(\mathcal{O}_{K}\right),$ it means that $e=2$ and $g=1.$ Since
 $\tau\left(2\right)=2,$ we obtain the equality in (2.2) and also in (2.1).\\
\textbf{Case 3.} If $p\mathcal{O}_{K}$$\notin$Spec$\left(\mathcal{O}_{K}\right)$ and $p\mathcal{O}_{K}\neq P^{2},$ with $P$$\in$Spec$\left(\mathcal{O}_{K}\right),$ it results that we have neither the possibility $e=g=1$ nor the possibility $e=2,$ $g=1.$ Applying these and the obvious inequality $$\tau\left(x\right)<x+1,\left(\forall\right) x\in\mathbb{N}^{*},$$ which implies $\tau\left(x\right)\leq x\left(\forall\right) x$$\in$$\mathbb{N}^{*},$ we are proving inequality (2.2) considering two subcases:\\
\textbf{Subcase 3.a)} If $e$$\in$$\mathbb{N}^{*},$ $e\geq 3$ and $g=1,$ inequality (2.2) is equivalent to $e\geq \tau\left(e\right),$ which is true. In fact, we have $e>\tau\left(e\right),$ for each $e$$\in$$\mathbb{N},$ $e\geq 3.$ Therefore, we obtain:\\
$\tau\left(p\mathcal{O}_{K}\right)>\tau^{\left(e\right)}\left(p\mathcal{O}_{K}\right) + \frac{\tau\left(p\mathcal{O}_{K}\right)}{\omega\left(p\mathcal{O}_{K}\right)}\cdot \frac{g}{e+1},$ $\left(\forall\right) e$$\in$$\mathbb{N}^{*},$ $e\geq 3$ and $g=1.$\\
\textbf{Subcase 3.b)} If $e;g$$\in$$\mathbb{N}^{*},$ $g\geq 2,$ we have: $0<\tau\left(e\right)<e+1.$ Since $g\geq 2,$ it results
$0<\left(\tau\left(e\right)\right)^{g-1}<\left(e+1\right)^{g-1}.$
Multiplying this inequality with inequality $0<\tau\left(e\right)\leq e,$ $\left(\forall\right) e$$\in$$\mathbb{N}^{*},$ we obtain 
$$\left(\tau\left(e\right)\right)^{g}< \left(e+1\right)^{g-1}\cdot e.$$
The above inequality is equivalent with
$$\tau\left(p\mathcal{O}_{K}\right)>\tau^{\left(e\right)}\left(p\mathcal{O}_{K}\right) + \frac{\tau\left(p\mathcal{O}_{K}\right)}{\omega\left(p\mathcal{O}_{K}\right)}\cdot \frac{g}{e+1}, \left(\forall\right) e,g\in\mathbb{N}^{*}, g\geq 2.$$
\end{proof}
When the extension of fields $\mathbb{Q} \subset K$ is a nontrivial Galois extension of degree odd, we are improving the Theorem 2.12 like this:

\begin{theorem}  
Let $n$$\in$$\mathbb{N}^{*}$ and let $K$ be an algebraic number field of degree odd $[K :
\mathbb{Q}] = n\geq3$ such that $\mathbb{Q} \subset K$ is a Galois extension. Let $p$ be a prime positive integer and let the ideal
$p\mathcal{O}_{K}=P^{e}_{1}\cdot P^{e}_{2}\cdot...\cdot P^{e}_{g},$ where $g$$\in$$\mathbb{N}^{*},$ $e$$\in$$\mathbb{N}^{*},$ and $P_{1},$ $P_{2},$..., $P_{g}$ are prime different ideals in the ring of integers $\mathcal{O}_{K}.$ If the ideal $p\mathcal{O}_{K}$ is ramified, then, we have the following inequality:
\begin{equation}{\label{3.3}}
\frac{e}{e+1}\geq \frac{\tau^{\left(e\right)}\left(p\mathcal{O}_{K}\right)}{\tau\left(p\mathcal{O}_{K}\right)} + \frac{1}{\omega\left(p\mathcal{O}_{K}\right)}\cdot \frac{g}{e+1}.
\end{equation}
\end{theorem}
\begin{proof}
Let $f$ be the the residual degree of $p.$ Since $\mathbb{Q} \subset K$ is a Galois extension, we have $f=f_{i}=\left[\mathcal{O}_{K}/P_{i} : \mathbb{Z}/p\mathbb{Z}\right],$ $\left(\forall\right)$ $i = \overline{1,g}.$ Since the ideal $p\mathcal{O}_{K}$ is ramified, it results $e>1.$ But $efg=n\geq 3$ and $n$ is an odd number, it results that $e\geq3.$\\
We have: $\tau\left(p\mathcal{O}_{K}\right)=\left(e+1\right)^{g},$ $\tau^{\left(e\right)}\left(p\mathcal{O}_{K}\right)=\left(\tau\left(e\right)\right)^{g},$ $\omega\left(p\mathcal{O}_{K}\right)=g.$ Therefore, the inequality (2.3) is equivalent with 

$$\frac{e}{e+1}\cdot\left(e+1\right)^{g}\geq \left(\tau\left(e\right)\right)^{g}+\frac{\left(e+1\right)^{g}}{g}\cdot \frac{g}{e+1} \Leftrightarrow$$
$$\left(e+1\right)^{g-1}\cdot \left(e-1\right)\geq \left(\tau\left(e\right)\right)^{g}.$$
In order to prove this last inequality, we are using the fact that $\tau\left(e\right)\leq e-1,$ $\left(\forall\right) e$$\in$$\mathbb{N}^{*},$ $e\geq3.$ Since $g$$\in$$\mathbb{N}^{*},$ it results that $\left(\tau\left(e\right)\right)^{g}\leq \left(e-1\right)^{g}.$ But 
$$\left(e-1\right)^{g}\leq \left(e+1\right)^{g-1}\cdot \left(e-1\right), \left(\forall\right) g\in\mathbb{N}^{*}.$$ 
So, by transitivity, we obtain:
$$\left(\tau\left(e\right)\right)^{g}\leq \left(e+1\right)^{g-1}\cdot \left(e-1\right), \left(\forall\right) e,g \in \mathbb{N}^{*}, e\geq3.$$
This last inequality is equivalent with inequality (2.3).\\
We remark that when $e=3, g=1$ we obtain the equality in (2.3).

\end{proof}
Now, we give some examples of fields extensions and of ideals which verify Theorem 2.13:\\
\smallskip\\
1) Let $\epsilon$ be a primitive root of order $3$ of the unity, let $p$ be a prime positive integer, $p$$\equiv$$2$ (mod $3$) and let $l$ be an integer such that $p|l.$ Let the cyclotomic field $K=\mathbb{Q}\left(\epsilon\right)$ and let the Kummer field $L=\mathbb{Q}\left(\epsilon; \sqrt[3]{l}\right).$ The fields extension $K\subset L$ is a Galois extension of degree $3.$ Since $p$$\equiv$$2$ (mod $3$), applying Theorem 1.4, it results that $p\mathbb{Z}\left[\epsilon\right]$ is a prime ideal. Since $p|l,$ it results that
the cubic character $\left( \frac{l}{p\mathbb{Z}\left[\epsilon\right]}%
\right) _{3}=0$ and according to Theorem 1.6, we obtain that the ideal $p\mathcal{O}_{L}=P^{3}_{1},$ where $P_{1}$$\in$Spec($\mathcal{O}_{L}$). So, the ideal $p\mathcal{O}_{L}$ is totally ramified in $\mathcal{O}_{L}$.
Using the notations from Theorem 2.13, we have: $g=1,$ $e=3,$ $\tau\left(p\mathcal{O}_{L}\right)=\left(e+1\right)^{g}=4,$ $\tau^{\left(e\right)}\left(p\mathcal{O}_{L}\right)=\left(\tau\left(e\right)\right)^{g}$$=2,$ $\omega\left(p\mathcal{O}_{L}\right)=g=1.$ It results
$\frac{e}{e+1}=\frac{3}{4}$ and $\frac{\tau^{\left(e\right)}\left(p\mathcal{O}_{L}\right)}{\tau\left(p\mathcal{O}_{L}\right)} + \frac{1}{\omega\left(p\mathcal{O}_{L}\right)}\cdot \frac{g}{e+1}=\frac{3}{4}.$ So, we have equality in (2.3).\\
2) Let $\xi$ be a primitive root of order $5$ of the unity, let $p$ be a prime positive integer, $p$$\equiv$$3$ (mod $5$) and let $l$ be an integer such that $p|l.$ Let the cyclotomic field $K=\mathbb{Q}\left(\xi\right)$ and let the Kummer field $L=\mathbb{Q}\left(\xi; \sqrt[5]{l}\right).$ The fields extension $K\subset L$ is a Galois extension of degree $5.$ Since $p$$\equiv$$3$ (mod $5$), and applying Theorem 1.4, it results that $p\mathbb{Z}\left[\xi\right]$ is a prime ideal. Since $p|l,$ it results that
the the $5$-power character $\left( \frac{l}{p\mathbb{Z}\left[\xi\right]}%
\right) _{5}=0$ and according to Theorem 1.6, we obtain that the ideal $p\mathcal{O}_{L}=P^{5}_{1},$ where $P_{1}$$\in$Spec($\mathcal{O}_{L}$). So, the ideal $p\mathcal{O}_{L}$ is ramified totally in $\mathcal{O}_{L}$.
Using the notations from Theorem 2.13, we have: $g=1,$ $e=5,$ $\tau\left(p\mathcal{O}_{L}\right)=\left(e+1\right)^{g}=6,$ $\tau^{\left(e\right)}\left(p\mathcal{O}_{L}\right)=\left(\tau\left(e\right)\right)^{g}$$=2,$ $\omega\left(p\mathcal{O}_{L}\right)=g=1.$ It results
$\frac{e}{e+1}=\frac{5}{6}$ and $\frac{\tau^{\left(e\right)}\left(p\mathcal{O}_{L}\right)}{\tau\left(p\mathcal{O}_{L}\right)} + \frac{1}{\omega\left(p\mathcal{O}_{L}\right)}\cdot \frac{g}{e+1}=\frac{1}{2}.$ So, we have
$$\frac{e}{e+1} > \frac{\tau^{\left(e\right)}\left(p\mathcal{O}_{L}\right)}{\tau\left(p\mathcal{O}_{L}\right)} + \frac{1}{\omega\left(p\mathcal{O}_{L}\right)}\cdot \frac{g}{e+1}.$$
\smallskip\\
The same two examples verify Proposition 2.9.

\end{document}